\documentclass[12pt,reqno,a4paper,twoside]{article}

\usepackage{amsmath,amsthm,amstext,amscd,amssymb,euscript,mathrsfs}
\usepackage{calrsfs}
\usepackage{epsfig}

\newcommand{\R}{\mathbb R}
\newcommand{\N}{\mathbb N}

\newcommand{\alphaa}{\kappa}

\renewcommand{\phi}{\varphi}

\newcommand{\pee}{\ensuremath{\mathbb{P}}}

\def\1{{\mathchoice {\rm 1\mskip-4mu l} {\rm 1\mskip-4mu l}
{\rm 1\mskip-4.5mu l} {\rm 1\mskip-5mu l}}}

\newtheorem{theorem}{{\small T}{\scriptsize HEOREM}}[section]
\newtheorem{corollary}{{\bf{\small C}{\scriptsize OROLLARY}}}[section]
\newtheorem{proposition}{{\bf{\small P}{\scriptsize ROPOSITION}}}[section]
\newtheorem{lemma}{{\bf{\small L}{\scriptsize EMMA}}}[section]
\newtheorem{remark}{{\bf{\small R}{\scriptsize EMARK}}}[section]
\newtheorem{definition}{{\bf{\small D}{\scriptsize EFINITION}}}[section]

\renewenvironment{proof}[1]
{\noindent{{\bf{\small{ P}{\scriptsize ROOF}}}.}\hspace{0.1cm} #1} {$\;\qed$\newline}

\newcommand{\beq}{\begin{eqnarray}}
\newcommand{\eeq}{\end{eqnarray}}

\newcommand{\ba}{\begin{align*}}
\newcommand{\ea}{\end{align*}}

\newcommand{\be}{\begin{equation}}
\newcommand{\ee}{\end{equation}}

\newcommand{\bl}{\begin{lemma}}
\newcommand{\el}{\end{lemma}}

\newcommand{\br}{\begin{remark}}
\newcommand{\er}{\end{remark}}

\newcommand{\bt}{\begin{theorem}}
\newcommand{\et}{\end{theorem}}

\newcommand{\bd}{\begin{definition}}
\newcommand{\ed}{\end{definition}}

\newcommand{\bp}{\begin{proposition}}
\newcommand{\ep}{\end{proposition}}

\newcommand{\bc}{\begin{corollary}}
\newcommand{\ec}{\end{corollary}}

\newcommand{\bpr}{\begin{proof}}
\newcommand{\epr}{\end{proof}}

\newcommand{\bi}{\begin{itemize}}
\newcommand{\ei}{\end{itemize}}

\newcommand{\ben}{\begin{enumerate}}
\newcommand{\een}{\end{enumerate}}

\newcommand{\caB}{{\mathcal B}}

\newcommand{\caE}{{\mathrsfs E}}

\newcommand{\caH}{{\mathcal H}}
\newcommand{\caI}{{\mathcal I}}

\newcommand{\caL}{{\mathcal L}}

\newcommand{\caS}{{\mathcal S}}

\begin{document}
\title{Gibbs-non-Gibbs transitions via large deviations:
computable examples}
\author{
Frank Redig$^{\textup{{\tiny(a)}}}$, Feijia Wang$^{\textup{{\tiny(b)}}}$\\
{\small $^{\textup{(a)}}$
Delft Institute of Applied Mathematics,}\\
{\small Technische Universiteit Delft}\\
{\small Mekelweg 4, 2628 CD Delft, Nederland}\\
{\small $^{\textup{(b)}}$ Mathematisch Instituut Universiteit Leiden}\\
{\small Niels Bohrweg 1, 2333 CA Leiden, The Netherlands}}

\maketitle

\begin{abstract}
We give new and explicitly computable examples of
Gibbs-non-Gibbs transitions of mean-field type, using
the large deviation approach introduced in \cite{efhr2}. These examples include Brownian motion with
small variance and related diffusion processes, such as the Ornstein-Uhlenbeck process, as well as birth and death processes.
We show for a large class of initial measures and diffusive dynamics both short-time conservation of Gibbsianness
and dynamical Gibbs-non-Gibbs transitions.
\bigskip
\noindent
\end{abstract}

\section{Introduction}
Starting from \cite{efhr} dynamical Gibbs-non-Gibbs transitions have been considered by several
authors, see e.g.\ \cite{ro}, \cite{kul}, \cite{rus}. In these studies,
one considers lattice spin systems started from a Gibbs measure $\mu$ at time zero and evolves it according to a
Markovian dynamics (e.g. Glauber dynamics)
with stationary Gibbs measure $\nu\not=\mu$. The question is then whether $\mu_t$, the time-evolved measure at
time $t>0$ is a Gibbs measure. Typically this is the case for short times, whereas for longer times,
there can be transitions from Gibbs to non-Gibbs (loss) and back from non-Gibbs to Gibbs (recovery).
The notion of a ``bad configuration'', i.e., a point of essential discontinuity of
the conditional probabilities of the measure $\mu_t$ is crucial here. Such a configuration $\eta_{spec}$
is typically identified by looking at the joint distribution of the system at time $0$ and at time $t$.
If conditioned on $\eta_{spec}$ the system at time zero has a phase transition, then typically
$\eta_{spec}$ is a bad configuration.

In the context of mean-field models, the authors in \cite{kuel} started with an analysis of the most
probable trajectories (in the sense of large deviations)
of a system conditioned to arrive at time $T$ at a given configuration.
The setting of \cite{kuel} is the Curie-Weiss model
subjected to a spin-flip dynamics. A Gibbs-non-Gibbs transition is in this context rephrased as a phenomenon
of ``competing histories'', i.e., for special terminal conditions $x_{spec}$ and
times $T$ not too small, multiple trajectories can minimize the rate
function, and these trajectories can be selected by suitably approximating
$x_{spec}$. These special conditionings leading to multiple histories
are the analogue of ``bad configurations'' (essential points of discontinuity
of conditional probabilities of the measure at time $t$) in the Gibbs-non-Gibbs transition
scenario.
This ``trajectory large deviation approach'' has then been studied in more generality, including
the lattice case,
in \cite{efhr2}.

In this paper, we apply the trajectory large deviation approach
in several examples, both for diffusion processes and for birth and death processes.
This leads to new and explicitly computable Gibbs-non-Gibbs transitions
of mean field type. For processes of diffusion type, we first treat an explicit
example for the rate function of the initial measure, and as dynamics Brownian motion
with small variance or Ornstein-Uhlenbeck process.
In all cases, we obtain the explicit form of the conditioned trajectories, and explicit
formulas for the bad configuration and the time at which it becomes bad.
In the case of general Markovian diffusion processes in a symmetric potential landscape,
we show under reasonable conditions short-time Gibbsianness as well as appearance
of bad configurations at large times.
Next, we treat the case of continuous-time random walk with small increments, as arises e.g.\
naturally in the context of (properly rescaled) population dynamics. In that case, the Euler-Lagrange
trajectories can be explicitly computed for some particular choices of the ``birth and death'' rates.
Constant birth and death rates are the analogue of the Brownian motion case, whereas linear birth
and death rates are the analogue of the Ornstein-Uhlenbeck process, but in that case the cost of
optimal trajectories becomes a much more complicated expression.

Our paper is organized as follows. In section \ref{fksec} we introduce some elements
or the Feng-Kurtz formalism, and define the notion of bad configurations in the present setting.
In section \ref{rcase} we treat diffusion processes with small variance, with an explicit
form for the initial rate function.
In section \ref{oe} we treat the case of Brownian motion dynamics with different cases for the
rate function of the initial measure.
Finally, in section \ref{bisec}, we treat one-dimensional random
walks with small increments, such as rescaled birth and death processes.

\section{The Feng-Kurtz scheme, Euler-Lagrange trajectories, bad configurations}\label{fksec}
We study Markov processes $\{X^n_t: 0\leq t\leq T\}$ taking values in $\R^d$, parametrized by
a natural number $n$.
This parameter tunes the ``amount of noise'' in the process, i.e.,
as $n\to\infty$, the process becomes deterministic, and the measure
on trajectories satisfies the large deviation principle with rate $n$
and with a rate function of the form
\be\label{lagrarate}
\caI(\gamma)= \int_0^T L(\gamma_s,\dot{\gamma}_s)ds
\ee
This means more precisely that
\be\label{largedevpath}
\pee\left(\{X^n_t:0\leq t\leq T\}\approx\gamma\right)\approx
\exp \left(-n \caI(\gamma)\right)
\ee
to be interpreted in the usual sense of the large deviation
principle with a suitable topology on the set of trajectories.
The form
\eqref{lagrarate} naturally follows from the Markov property.

Notice that the form of the rate function does not depend on the
choice of this topology. So one usually starts with the weakest
topology, i.e., the product topology, and then, if possible,
strengthens the topology by showing exponential tightness. See
\cite{dz} for an illustration of this strategy in the context of
theorems like Mogulskii's theorem.

Since in this paper we are only interested in finding out
optimal trajectories, i.e., minimizers of the rate function
over a set of trajectories with prescribed terminal condition
and open-end condition, we will not have to worry about
the strongest topology in which the large deviation principle
\eqref{largedevpath} holds, but we are rather after (as explicit
as possible)
solutions of Euler-Lagrange problems associated to the rate
function.

In \cite{fk} a scheme is given to compute
the ``Lagrangian'' $L$, see also \cite{efhr2} for an illustration
of this scheme in the large-deviation view on Gibbs-non-Gibbs transitions. First one computes the ``Hamiltonian''
\be\label{hamfeng}
\caH (p,x) = \lim_{n\to\infty}\frac1n e^{-n\langle p,x\rangle}\caL_n e^{n\langle p,x\rangle}
\ee
where $\caL_n$ is the generator of the process
$\{X^n_t: 0\leq t\leq T\}$ (working on the $x$-variable, where
$p\in\R^d$ is the ``momentum'' and where $<.,.>$ denotes inner product. .
Under regularity conditions on $\caH (p,x)$ (e.g.\ strict convexity), the associated Lagrangian is then
given by the Legendre transform
\be\label{legenhamlag}
L(x,v)=\sup_{p\in\R^d}\left(<v,p>- \caH(x,p)\right),
\ee
As an example, consider
\[
X^n_t= n^{-1/2}B_t
\]
with generator
\[
\caL_n= \frac1{2n} \Delta
\]
then we have
\[
\caH(x,p)= \frac{p^2}{2}
\]
and associated Lagrangian
\[
L(x,v)=\frac{v^2}{2}
\]
which produces the rate function of the well-known Schilder's theorem
\[
\pee\left(\{X^n_t:0\leq t\leq T\}\approx \gamma\right)
\approx\exp\left(-\frac{n}{2}\int_0^T \dot{\gamma}_s^2 ds\right)
\]
To proceed, we also want the initial point of our process
to have some fluctuations.
More precisely, we need for the starting point of our process
an initial measure
$\mu_n$ (depending on $n$) on $\R^d$, satisfying the large deviation principle
with rate $n$ and rate function $i(x)$, i.e., in the sense of
large deviations,
we assume
\be\label{initialra}
\pee\left(X^n_0\in A\right)=\mu_n(A)\approx \exp(-n\inf_{x\in A} i(x))
\ee
We call the triple $(\{X^n_t:0\leq t\leq T\}, L,i)$
a stochastic system with small noise.

We continue now with the definition of a bad configuration
in this framework. This is motivated by the definition
of a bad configuration in the context of mean-field models
\cite{kuel}, and can be viewed as the large-deviation rephrasing
of ``a phase transition at time zero conditioned on a
special configuration at time $T$''.
\bd
Let $(\{X^n_t:0\leq t\leq T\},L,i)$ be a stochastic
system with small noise.  We say that a point $b\in\R^d$
is bad at time $T$ if the following two conditions hold.
\ben
\item Conditional on $X^n_T=b$, $X^n_0$
does not converge (as $n\to\infty$)
to a pointmass in distribution.
\item There exist two sequences $b^+_k\to b$, $b^-_k\to b$
and $\delta>0$
such that the variational distance between the distribution
$\mu(0,T;b^+_k)$ of $X^n_0|X^n_T=b^+_k$ and
the distribution
$\mu(0,T;b^-_k)$ of $X^n_0|X^n_T=b^-_k$ is
at least $\delta$ for $k$ large enough.
\een
\ed
The simplest example which follows also the most common scenario is
where the distribution of
$X^n_0|X^n_T=b$
converges to $\frac12(\delta_{-a}+\delta_{a})$ and for $c>b$
$X^n_0|X^n_T=c$
converges to $\delta_{\alpha(c)}$ where
$\alpha(c)\to a$ as $c\downarrow b$, whereas
for $c<b$
$X^n_0|X^n_T=c$
converges to $\delta_{\alpha'(c)}$ where
$\alpha'(c)\to -a$ as $c\uparrow b$.
This means that conditioned to be at time $T$ at location $b$,
the process has two ``favourite'' intial spots, which can be
``selected'' by approaching $b$ from the right or from the left.

This is the analogue of a phase transition, where the phases
can be selected by appropriately approximating the bad configuration, see
\cite{efhr}.

\section{Diffusion processes with small variance conditioned on the future}\label{rcase}
In this section we present examples where
$X^n_t$ is a diffusion process.
We show also how from the large deviation approach we gain a new
understanding of ``short-time Gibbsianness''
for a general class of drifts of the diffusion, or
initial rate functions.

\subsection{Brownian motion}
\begin{figure} 
\centering
\includegraphics[height=80mm]{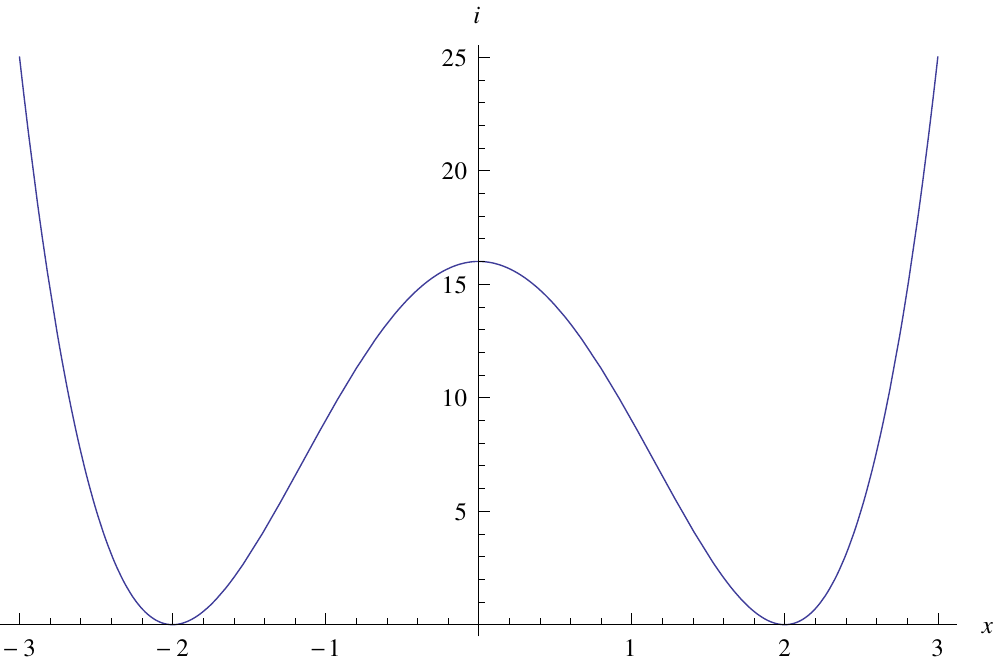}
\caption{$i(x)= (x^2-a^2)^2$ with $a=2$}
\label{ratefunc1}
\end{figure}

To start with, we consider Brownian motion with small variance
$\frac1n$ starting
from an initial distribution satisfying the large deviation principle
(with rate $n$)
with a non-convex rate function having two mimina at locations $-a,a$,
with $a>0$.
More precisely, we consider the process
\be\label{brown}
X^n_t =\frac{1}{\sqrt{n}} B_t
\ee
starting from an initial distribution $\mu_n$ such that, informally written,
\be\label{initial}
\pee\left(X^n_0\in dx\right)=\mu_n (dx) \approx e^{-n i(x)} dx
\ee
For $i$ we make the explicit choice:
\be\label{rate}
i(x)= (x^2-a^2)^2
\ee
i.e., a non-convex function non-negative with zeros at $-a,a$ and maximum at $x=0$ ($i(x)$ with $a=2$ is plotted in Figure \ref{ratefunc1}).

This specific choice is for the sake of {\em explicit analytic computability}
but many results are true for a general class of rate functions
that have a similar graph with two zeros located at $-a,a$ and
a maximum at zero.

More formally, we require that the sequence of initial probability  measures $\{\mu_n,n\in\N\}$
satisfies the large deviation principle with rate function $i$ given by
\eqref{rate}. Such rate functions arise naturally in the context of mean-field
models with continuous spins and spin-Hamiltonian depending on the magnetization.

We are then interested in the most probable trajectory $\gamma$
with initial
point distributed according to $\mu_n$ and
final point $\gamma_T=0$.
More precisely, by application of Schilder's theorem, the trajectory
$\{X^n_t: 0\leq t\leq T\}$ satisfies the LDP with rate function
\be\label{pad}
\caI (\gamma) =\frac12 \int_0^T \dot{\gamma}_s^2 ds + i(\gamma_0)
\ee
The optimal trajectory we are looking for is hence
\[
\arg\min \{\caI (\gamma) : \gamma(T)= 0\}
\]
The Euler-Lagrange trajectories (extrema of the cost
$\frac12\int_0^T\dot{\gamma}_s^2 ds$ corresponding to $\caI(\gamma)$) are linear in $t$:
\[
\gamma_t = A+ Bt
\]
By the terminal condition $B= -A/T$.

The cost $\caI (\gamma)$ of this trajectory can then be rewritten
as a function of the starting point $\gamma_0=A$:
\be\label{kotso}
\caE_{0,T}(A):=\caI (\gamma)= A^4- 2a^2 A^2 + a^4 + \frac12 (-A/T)^2 T 
= A^4 +\alpha(a,T) A^2 + a^4
\ee
with
\be\label{alp}
\alpha(a,T)= \left(\frac1{2T}-2a^2\right)
\ee
The behavior of this cost depends on
the sign of $\alpha$.
If $\alpha \geq 0$, then there is a unique minimum at $A=0$, this
case corresponds
to
\[
T \leq \frac{1}{4a^2}:= T_{crit}
\]
If $\alpha <0$ then there are two mimima $A=A_{\pm}$ given by
\be\label{biba}
A_{\pm} = \pm \sqrt{-\alpha(a,T)/2}=\pm \sqrt{a^2-(4T)^{-1}}
\ee

We thus conclude that,
as $n\to\infty$, the starting point is most probably $0$ for small
$T$ and most (and equally) probably $A_{\pm}$ for large $T$, which
converges to $\pm a$ when $T\to\infty$.
Hence we have non-uniqueness of histories.

Let us denote $\mu(n,T,0)$ the distribution
of $X^n_0$ conditioned on $X^n_T=0$.
Then we have
\ben
\item Small times, unique history. If $T \leq T_{crit}$ then
\[
\lim_{n\to\infty}\mu(n,T,0)=\delta_0.
\]
\item Large times, non-unique history. If $T > T_{crit}$ then
\[
\lim_{n\to\infty}\mu(n,T,0)= \frac12(\delta_{A^+}+\delta_{A^-}).
\]
\item Limit of large times
\[
\lim_{T\to\infty}\lim_{n\to\infty}\mu(n,T,0)\to \frac12(\delta_{a}+\delta_{-a}).
\]
\een

Let us now condition on $X^n_T=b\not=0$.
Then the most probable trajectory is still a straight line
$\gamma^b_t=A+Bt$
but now with terminal condition $A+BT=b$, i.e.,
$B= (b-A)/T$. It has cost expressed in terms of the starting point $\gamma_0=A$
\be\label{kotsb}
\caE_{b,T} (A)= A^4 +\alpha(a,T) A^2 + a^4-\frac{b}{T}A+ \frac{b^2}{2T}
\ee
This is the cost function $\caE_{0,T}(A)$ of \eqref{kotso} plus a linear term
$-\frac{b}{T}A +\frac{b^2}{2T} $.
Minimization of $\caE_{b,T}(A)$ leads to the equation
\be\label{boel}
4A^3 +2\alpha A = \frac{b}{T}
\ee

We then have two cases:
\ben
\item $\alpha \geq 0$, i.e., $T\leq T_{crit}$. Equation \eqref{boel}
has a unique real solution, corresponding to a unique minimum
$A_b$
of $E_b(A)$. This minimum converges to zero as $b\to 0$. Hence,
$0$ is good for $T\leq T_{crit}$.
\item $\alpha <0$. Equation \eqref{boel} has three real
solutions. For $b>0$ we have one positive and two negative solutions.
The positive solution denoted $A(+,b,T)>\sqrt{-\alpha/2}$ gives the minimum. The negative solutions correspond
to a maximum and a local minimum. For $b<0$ the situation
exactly the opposite: the unique negative solution $A(-,b,T)<-\sqrt{-\alpha/2}$
correspond to the global minimum whereas the two positive solutions
give a maximum and a local minimum. Hence $0$ is bad
for all $T>T_{crit}$
\een
In particular, for the $T\to\infty$ the positive, resp.\ negative
minimum of the rate function of the distribution at time zero is
selected by taking the right or left limit of the conditioning.
\[
\lim_{c\to 0 , b>0}\lim_{T\to\infty}\pee (X^n_0 = \cdot| X^n_T = c)=\delta_a
\]
and, similarly
\[
\lim_{c\to 0 , b<0}\lim_{T\to\infty}\pee (X^n_0 = \cdot| X^n_T = c)=\delta_{-a}
\]

Summarizing our findings, let us denote $\caB_T$ the set of bad configurations
then we have
\bt
\ben
\item Short times: no bad configurations.

For $T \leq \frac{1}{4a^2}$, $\caB_T=\emptyset$.
\item Large times: unique bad configuration.
For $T > \frac{1}{4a^2}$, $\caB_T=\{0\}$
\een
\et

\subsection{Brownian motion with constant drift}
The case of Brownian motion with constant drift $V>0$ is treated similarly.
The Euler-Lagrange trajectories are once more linear in $t$,
but the cost is now
\[
i(\gamma_0) + \frac{1}{2}\int_0^T(\dot{\gamma}_s- V)^2 ds
\]
which for $\gamma_t= A+Bt$ ending in $\gamma_T= b$ can
be computed explicitly and gives
\[
E_{b,V} (A)= \frac12\left(\frac{b-A}{T}- V\right)^2 + i(A)
\]
of which a similar analysis can be given. In particular,
choosing $b=VT$ we see that $VT$ is the cost is
identical to the zero drift case conditioning to be at zero
at time $T$, and hence this is a bad point for $T> T_{crit}$, where
$T_{crit}$ is the same critical time as for the zero drift case.
The analysis around this bad point is identical.
Notice that the ``limiting deterministic dynamics'' is
$\dot{x}= V$ and the bad point $x_{spec}= VT$ is precisely where this
dynamics ends up at time $T$ when started from zero.
\subsection{Other rate functions for the initial measure and corresponding behavior of Brownian motion}\label{oe}
We now consider other possible scenarios for different rate functions associated to the initial measure,
and for the Brownian motion with small
variance as dynamics. The starting measure $\mu_n(dx)$ satisfies the large deviation principle
with rate function $i(x)$.
As a consequence, the minimizing trajectory to arrive at position $b$ at time $T$
is $\gamma_t = Bt+A$ with $B=(b-A)/T$ and has cost
\begin{equation}
\label{ds}
\caE_{b,T}(A)= \frac {(b-A)^2}{2T} + i(A)
\end{equation}
The following scenarios can then occur
\ben
\item {\bf $i(A)$ is strictly convex: no bad configurations.} Indeed, in that case $\caE(A)$ is also
strictly convex (as a sum of two strict convex function) and hence has a unique minimum.
In this scenario, there are no bad configurations, and the optimal conditioned trajectory
is always unique. This corresponds to ``high temperature initial measure''
and ``infinite-temperature dynamics'', which always conserves Gibbsianness.

\begin{figure} 
\centering
\includegraphics[height=80mm]{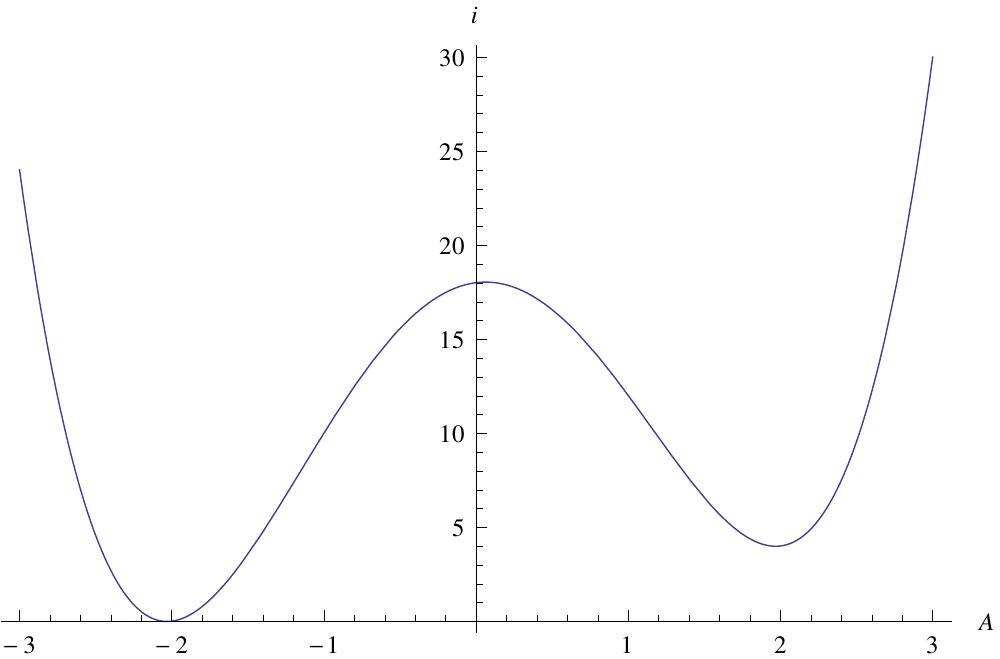}
\caption{$i(A)=(A^2-a^2)^2+A+r$ with $a=2$ and $r=2.01539$}
\label{ratefunc2}

\centering
\includegraphics[height=20mm]{blank.pdf}

\centering
\includegraphics[height=80mm]{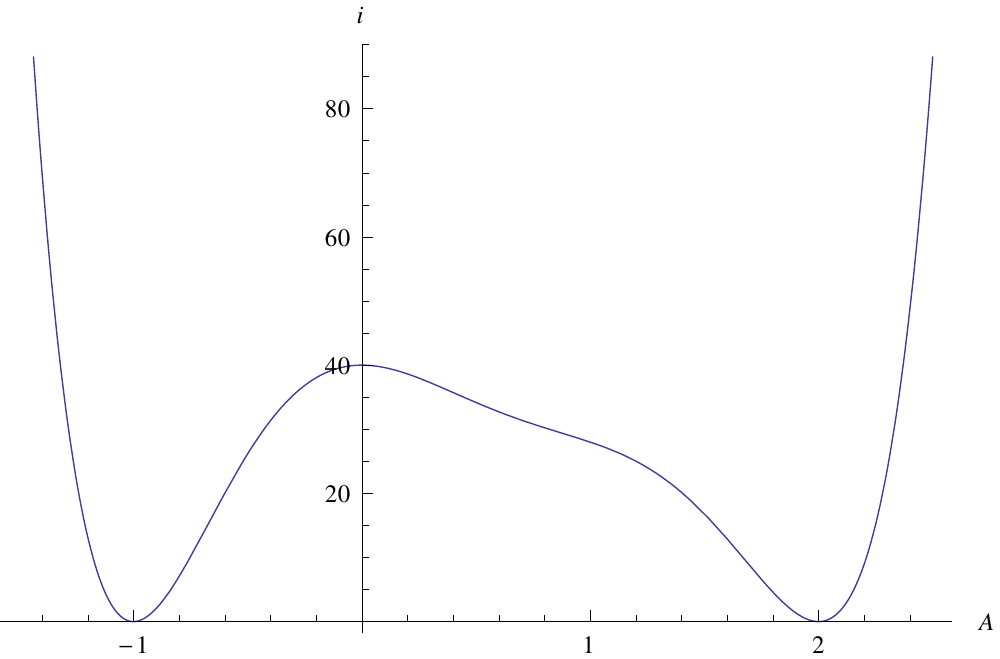}
\caption{$i(A)=7A^6-24A^5+9A^4+38A^3-42A^2+40$}
\label{ratefunc3}
\end{figure}

\item {\bf Initial field: loss without recovery, with a ``compensating'' bad configuration.}
As an example we can take $i(A)=(A^2-a^2)^2+A+r$. For $a>1$, this rate function has one local minimum in the vicinity of $x=a$, a maximum in the vicinity of $x=0$ and its
(absolute) minimum in the vicinity of $x=-a$. This corresponds to an initial field (favorizing the minimizer $x=-a$). $i(A)$ with $a=2$
 with $r=2.01539$ is plotted in Figure \ref{ratefunc2}.
The minimization of $\caE_{b,T}(A)$
leads to the equation
\begin{equation} \label{field}
4A^3+2\alpha(a,T) A = \frac{b}{T} - 1
\end{equation}
By an analysis of \eqref{field} similar for \eqref{boel}, we obtain that
there is no bad point when $T\leq T_{crit}=\frac{1}{4a^2}$, but $b=T$ is bad for all $T>T_{crit}$. The bad point ``compensates''
the initial field, and therefore has to become larger (and positive) when time $T$ increases.

\item {\bf Non-symmetric rate function.}
To see that the symmetry of the initial rate function is not a necessary requirement to produce bad configurations, we have
the following example.
Let $i(A)= 7A^6-24A^5+9A^4+38A^3-42A^2+40$ (see Figure \ref{ratefunc3}). This rate function has two global minima at $x=-1$ and
$x=2$ and one maximum at $x=0$. The cost function corresponding to trajectories arriving
at $b$ at time $T$ is
\[
\caE_{b,T} (A)= 7A^6 -24A^5 + 9A^4 + 38 A^3 + \left(\frac1{2T}-42\right)A^2 -\frac{b}{T} A + \frac{b^2}{2T} + 40
\]
For fixed $b$, and $T$ large enough, this function has two local minima, located
at $A^1(b,T)<A^2(b,T)$. Let us denote, for fixed $T$,
\[
D_T (b)= \caE_{b,T} (A^1(b,T))-\caE_{b,T} (A^2(b,T))
\]
If as a function of $b$, $D_T$ changes sign, by continuity, there must
be a value of $b^*$ where $D_T (b^*)=0$, i.e., where the minima of $\caE_{b^*,T}$
are at equal height. This $b^*$ is then a bad point at time $T$.
For $T=1$ we have $D_T(0.499) \approx -0.00182497<0$ and $D_T (0.4999)\approx 0.000868034>0$, so
at $T=1$, there is a bad point at $b^*\in (0.499, 0.4999)$. We observe that $b^*$ is $T$ dependent and tends to $0.5$ as $T$ increases. From numerical
oberservations, we have $b^*\in (0.4999, 0.49999)$ for $T=4$, $b^*\in (0.49999, 0.499999)$ for $T=39$ and $b^*\in (0.499999, 0.4999999)$ for $T=1000$.

\item {\bf General symmetric rate function.} For {\em any rate function $i(A)$ which is symmetric with respect to $x=0$} and which has minima
for $A\neq 0$, $b=0$ is bad when $T$ is large enough. Indeed, the cost to arrive at $0$ is from \eqref{ds}: $i(A)+\frac{A^2}{2T}$ which has a non zero minimum
as soon as $T$ is large enough.

\item {\bf General short-time Gibbsianness.} For {\em every rate function $i$ which is twice differentiable and
its second derivative is continuous and bounded from below}, we show that for $T$ small
enough there is a unique minimum $A_b$ of $\caE_{b,T}(A)$. This
is the analogue of ``short-time'' Gibbsianness obtained
in the lattice case via cluster expansions
\cite{leny}
or conditional Dobrushin uniqueness
\cite{opokuk} and can be
proved as follows.

We look at the equation (see \eqref{ds})
\begin{equation} \label{min-eq}
i'(A)=-\frac{A}{T}+\frac{b}{T}=:f(A).
\end{equation}

Put $d=\inf_{A}i''(A)$. Then we conclude, for
\begin{equation}
T<-\frac{1}{d},
\end{equation}
that \eqref{min-eq} has only one real solution $A_b$. Indeed, look at any two adjacent intersection points $A_1$ and
$A_2$ of $i'(A)$ and $f(A)$ if there were more than one real
solution for \eqref{min-eq}. By the intermediate value
theorem, we get
\begin{equation}
\min(i''(A_1),i''(A_2))<-\frac{1}{T}<d=\inf_{A}i''(A).
\end{equation}
This is a contradiction. And further because $i''(A_b)>-\frac{1}{T}$, we
have
\begin{equation}
\caE_{b,T}''(A_b)=i''(A_b)+\frac{1}{T}>0.
\end{equation}
Therefore $A_b$ is a minimum.

\item {\bf Non-Gibbsianness for all times.} An example where $b=0$ is bad for all $T>0$ is $i(0)=0$, $i(A)=\int_0^{|A|}|x|\cos^2\frac{1}{x} dx$ for $A\neq0$.
This follows from the facts that $i''(A)$ is not bounded from below when $A\rightarrow 0$ and $i(A)$ is symmetric about $A=0$.
To see that indeed for all $T>0$ $0$ is a bad point, we see that the line $f(A)=-A/T$ always intersects the
graph of the derivative of the rate function.
\een

\section{Ornstein-Uhlenbeck process }
As a second example, we consider the process $X^n_t$ to be the solution of
\[ \label{OU}
dX_t= -\alphaa X_t + \frac{1}{\sqrt{n}} dB_t
\]
and the initial point distributed as in the previous
section, in \eqref{initial}, \eqref{rate}.

The cost function for the large deviation principle of the
trajectories now becomes
\be
\caI (\gamma)= i(\gamma_0) +\frac12 \int_0^T (\dot{\gamma}_s +\alphaa \gamma_s)^2 ds
\ee
The Euler-Lagrange trajectories extremizing
$\frac12 \int_0^T (\dot{\gamma}_s +\alphaa \gamma_s)^2 ds$
are given by
\[
\gamma_t= Ae^{\alphaa t} + Be^{-\alphaa t}
\]
by the terminal condition $\gamma_T=0$ we have
\[
\gamma_t = -Be^{-2\alphaa T}e^{\alphaa t} + B e^{-\alphaa t}
\]
the cost function for such a trajectory can then explicitly be evaluated and
gives
\be\label{boemba}
\caE_{0,T}(B)= c_1 B^4 + c_2 B^2 + c_3
\ee
where
\beq
c_1 &=& (1-e^{-2\alphaa t})^4
\nonumber\\
c_2 &=& \left(-2a^2 (1-e^{-2\alphaa t})^2 + \alphaa e^{-2\alphaa t}(1-e^{-2\alphaa t})\right)
\nonumber\\
c_3&=& a^4
\eeq
A similar analysis as in the previous section can now be started.
We have a unique minimum at $B=0$ of the cost function $\caE$ for
\be\label{orntime}
T\leq T_{crit}:=-\frac{1}{2\alphaa} \log\left(\frac{2a^2}{2a^2+\alphaa}\right)
\ee
and for $T>T_{crit}$, $0$ becomes the unique bad point for this process.

The cost of an optimal trajectory ending up at $b$ at time $T$ can also be expressed as a function of the starting
point $\gamma_0$, which gives the explicit expression
\be\label{inexp}
\caE_{b,T} (\gamma_0)= i(\gamma_0) + \frac{\kappa }{e^{2\kappa T}-1} (\gamma_0- b e^{\kappa T})^2
\ee

\begin{figure} 
\centering
\includegraphics[height=80mm]{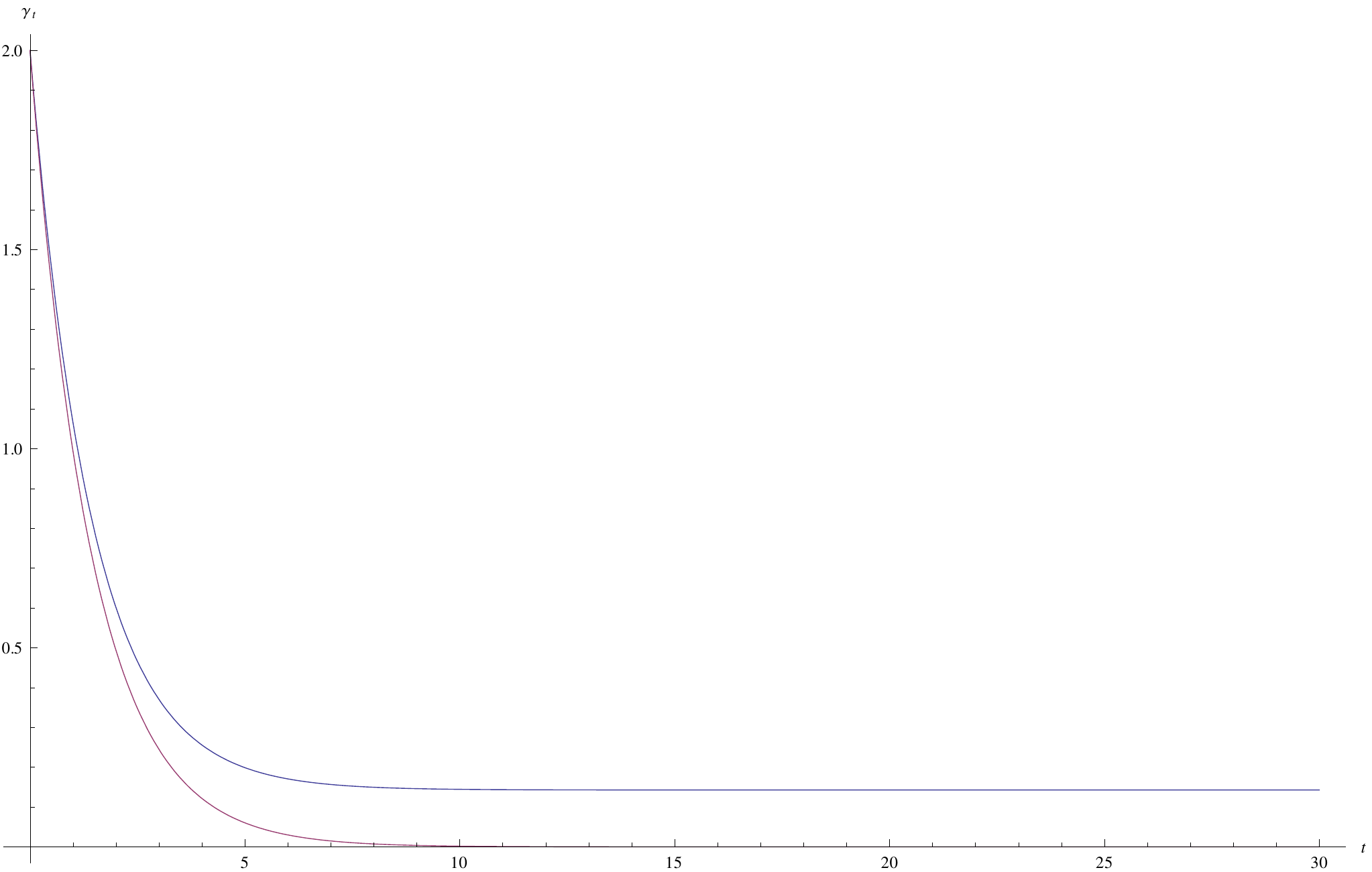}
\caption{A limiting process (the purple line) with $a=2,\kappa=0.7, T=30$, hence $\gamma_0^{+}\approx2.0$, and 
a corresponding conditioned process (the blue line) with $E=0.1$, $b\approx0.142857$}
\label{ou.pdf}
\end{figure}

\subsection{Ornstein-Uhlenbeck process with constant external field}
The equation for the process $X^n_T$ then reads
\be\label{orfie}
dX^n_t= (-\kappa X_t + E)dt + \frac{1}{\sqrt{n}} dB_t
\ee
where
$E>0$ is a constant representing a (constant) external field.
As rate function of the initial measure we choose as before
\eqref{rate}. The cost of the trajectory is now given by
$\int_0^T L (\gamma_s,\dot{\gamma}_s) ds$ with $L (\gamma_s,\dot{\gamma}_s) = (\dot{\gamma_s} +\kappa\gamma_s-E)^2$.
The Euler-Lagrange trajectories are of the form
\[
\gamma_t = Ae^{\kappa t} + Be^{-\kappa t} + \frac{E}{\kappa}
\]
The trajectory cost of an Euler-Lagrange trajectory is given by
$2A^2 (e^{2\kappa T}-1)$. From this, we derive that
the total cost of a trajectory to end up at time $T$ in $\gamma_T=b$ is given, as a function
of $\gamma_0$, by
\[
\caE^E_{b,T}(\gamma_0)= i(\gamma_0)+\frac{\kappa \left( \gamma_0 - \left(b-\frac{E}{\kappa}\right)e^{\kappa T}-\frac{E}{\kappa}\right)^2}{e^{2\kappa T}-1}
\]
The same analysis can then be performed. The ``critical'' time at which a unique bad point starts to appear is the same as in the
zero field case, i.e., given by
\eqref{orntime}. This bad point is given by
\be\label{badorn}
b= \frac{E}{\kappa}(1-e^{-\kappa T})
\ee
which corresponds to the point at which the deterministic evolution
$\dot{x}_t= -\kappa x_t +E$ arrives when starting from $x_0=0$.
Notice that total cost to arrive at this bad point $b$ is given by
\[
i(\gamma_0)+\frac{\kappa\gamma_0^2}{e^{2\kappa T}-1}
\]
which is symmetric around $\gamma_0=0$. Moreover, for $T$ large
the path cost contribution
which is equal to
$\frac{\kappa\gamma_0^2}{e^{2\kappa T}-1}$
vanishes exponentially fast, and hence for large $T$ two minima exist.

The corresponding optimal trajectories to arrive at the bad point $b$ are starting from
\[
\gamma_0^{\pm} =\pm \sqrt{ a^2 - \frac{\kappa}{2(e^{2\kappa T} -1)}}
\]
and explicitly given by
\begin{eqnarray*}
\gamma_t = \left(b-\frac{E}{\kappa}\right)\frac{\sinh (\kappa t)}{\sinh(\kappa T)}
+\left(\gamma_0-\frac{E}{\kappa}\right)\frac{\sinh (\kappa(T- t))}{\sinh(\kappa T)}
+\frac{E}{\kappa}
\end{eqnarray*}
The trajectory with plus
resp.\ minus sign can be selected by conditioning to arrive at $b^+>b$, resp.
$b^-< b$, and letting $b^+\to b$, resp.\ $b^-\to b$. Here we plot a limiting process with $a=2,\kappa=0.7, T=30$, hence $\gamma_0^{+}\approx2.0$, and 
a corresponding conditioned process with $E=0.1$, hence $b\approx0.142857$, see Figure \ref{ou.pdf}.

\subsection{General drift.}
Let us now consider the process $X^n_t$ with a general
drift $f(x)$ and variance $\frac1n$, i.e., the solution
of
\[
dX_t =-f(X_t) dt +\frac{1}{\sqrt{n}} dB_t
\]
We assume
$f:\R\to\R$ to be Lipschitz, and odd: $f(-x)=-f(x)$.
For the rate function of the initial point $X^n_0$ we choose
as before \eqref{initial}, \eqref{rate}.
The rate function of the trajectory is now given by
\be
\caI (\gamma)= \frac12\int_0^T (\dot{\gamma}_s + f(\gamma_s))^2 ds
\ee
and the minimization problem for the optimal trajectory
ending at zero $\gamma_T=0$ becomes now to find
\be\label{blurp}
\arg\min \{ \caI (\gamma) + i(\gamma_0): \gamma_T=0\}
\ee
The Euler-Lagrange equations for miminal cost trajectories are given by
\[
\frac{d^2\gamma_s}{ds^2}= f(\gamma_s)f' (\gamma_s)
\]
These equations correspond to classical motion in
a potential $U$ satisfying $U'=- f f'$, which gives as a possible
choice
$U= -\frac12 f^2$. Notice that this formal potential $U$ has no
physical meaning, but we need it if we want to translate the framework
of the Euler-Lagrange equations to Hamilton equations.
Indeed, the corresponding Hamiltonian is
\be\label{ham}
H(p,q)= \frac{p^2}{2} - U(q)
\ee
In particular, under the Euler-Lagrange equations,
\be\label{const}
\frac{\dot{\gamma}^2_t}{2}- \frac12 (f(\gamma_t))^2= E
\ee
is a constant of motion.
Further, we have the open-end and terminal condition
\beq
i'(\gamma_0)&=& \dot{\gamma_0} + f(\gamma_0)
\nonumber\\
\gamma_T&=&0
\eeq

We can think of these equations as having $\gamma_0$ and $E$ as
parameters. The terminal condition gives then a relation
between $E$ and $\gamma_0$.
Notice that the trajectory of zero-energy,  $E=0$, $\gamma\equiv 0$
is always a solution since $f(0)=0$.
We want to show that under some reasonable assumptions,
for
$T$ small, it is the only solution.
For this we make the following assumptions.
Call
$\caS_T(E)$ the collection of all trajectories $\gamma:[0,T]\to\R$ ending
at $0$, i.e., with $\gamma_T=0$ and with
``energy'' $E$, i.e., such that
\[
\frac{\dot{\gamma}^2_t}{2}- \frac12 (f(\gamma_t))^2= E
\]
for all $0\leq t\leq T$. We impose now the following conditions.
\begin{enumerate}\label{coco}
\item There exist a function
$\phi:\R\to [0,\infty)$ and $T_0>0$ and a constant $C>0$ such that
$\phi(0)=0$, $\phi(E)>0$ for $E\not=0$ such that for all $T\leq T_0$ and
for all $\gamma\in \caS_T(E)$, $\gamma_0\dot{\gamma}_0<0$,
\be\label{minaf}
|\dot{\gamma}_0|\geq \phi(E)
\ee
and
\be\label{maxt}
|\gamma_0|
\leq C\phi(E) T
\ee
\item The drift function $f$ is locally monotone around $0$, i.e.,
there exist $x_0$ such that $f$ restricted to $[0,x_0]$, $[-x_0,0]$ is monotone.
\end{enumerate}
The first condition states that if $T$ is small, and one wants
to end at $\gamma_T=0$ from $\gamma_0>0$, then the derivative at zero should be negative,
or vice versa. The second part of the condition states that there
exist lower bounds for the derivative and upper bounds for $\gamma_0$.

Coming back to the previous examples: for the Brownian motion case,
for all $\gamma\in \caS_T(E)$ we have
$\gamma_t= \pm\sqrt{2E} (t-T)$ hence and for
$\gamma\in \caS_T(E)$ we have $\gamma_0= \mp\sqrt{2E}T$,
$\dot{\gamma}_0=\pm\sqrt{2E}$, and
we can choose $\phi(E)= \sqrt{2E}$.
For the Ornstein-Uhlenbeck case
we have
$\gamma_t =B (e^{-\alphaa t}-e^{-\alphaa (2T-t)})$,
$E=-2AB\alpha$ and if $\gamma_T=0$ we
find
$\gamma_0= \sqrt{2E/\alphaa} \sinh(\alphaa T)$,
$\dot{\gamma}_0= -\sqrt{2E/\alphaa} \cosh(\alphaa T)$
which clearly satisfies the conditions, with the $\phi=\sqrt{2E/\alphaa}$.

The open-end condition requires
\[
\dot{\gamma}_0 + f(\gamma_0)= 4\gamma_0 (\gamma_0^2-a^2)
\]
Hence, for $\gamma\in \caS_T(E)$ such that $\gamma_0>0$:
\beq
-\phi(E)&\geq& \dot{\gamma}_0
\nonumber\\
&\geq&
4\gamma_0 (\gamma_0^2-a^2)-f(C\phi(E) T)
\nonumber\\
&\geq &
4C\phi(E) T (C^2\phi(E)^2 T^2-a^2)-f(C\phi(E) T)
\eeq
which is clearly a contradiction for $T$ sufficiently small.
Hence for $T$ sufficiently small, there do not exist $E\not=0$ with
$\gamma\in S_T (E)$. As a consequence, under these assumptions,
for small $T$ the zero trajectory is the only solution
of the minimization problem \eqref{blurp}.

For large times, if we assume that the drift is such that from any starting point one can travel to the origin at arbitrary small cost
if one has sufficient time, i.e., for all $x_0>0$,
\[
\lim_{T\to\infty} \inf \left\{ \int_0^T(\dot{ \gamma}_s + f(\gamma)_s)^2 ds: \gamma_0= x_0, \gamma_T=0\right\}=0
\]
then this implies that for $T$ large enough that there exists $x_0\not= 0$
and a trajectory $\gamma$ starting from $x_0$ such that $i(x_0)<i(0)/2$ and
\[
\left\{ \int_0^T(\dot{ \gamma}_s + f(\gamma)_s)^2 ds: \gamma_0= x_0, \gamma_T=0\right\} < i(0)/2
\]
this trajectory $\gamma$ clearly has lower cost than the zero trajectory, and
by symmetry, $-\gamma$ is a trajectory with identical cost. Therefore, $0$ becomes a bad point.


\section{Approximately deterministic walks in $d=1$}\label{bisec}
An ``approximately deterministic random walk''
is a continuous-time random walk with small increments performed
at high rate, i.e.,
a random walk $X^N_t$ on $\R$ that, starting at $X_0=x$ makes increments of size
$\pm 1/N$ with rates $Nb(x)$, resp.\ $Nd(x)$. In other words,
$X^N_t$ is a Markov process on $\R$ with generator
\be\label{birthgen}
\caL_N f(x)= Nb(x) \left(f\left(x+\frac1N\right)-f(x)\right)
+ N d(x)\left(f\left(x-\frac1N\right)-f(x)\right)
\ee

Such walks arise naturally in the context of population dynamics
see e.g. \cite{etheridge}. The notation $b(x)$ and $d(x)$
is also reminiscent of this interpretation and we will call these quantities birth resp.
death rates.

We ask then the same large deviation question, i.e.,
we start the process $X^N_t$ from an initial distribution
$\mu_N$ satisfying the large deviation principle with rate
function \eqref{rate} -or some natural modification of it if we
have to restrict the state space-
and look for the minimizing trajectory(ies)
that end at time $T$ at the origin (or at a more general bad point
if the dynamics has a drift see later).

The large deviation function for the trajectories can be computed
using the Feng-Kurtz scheme, i.e., denoting $f^N_p(x)= e^{Npx}$
we computes the Hamiltonian
\be\label{hamil}
\caH(x,p)=\lim_{N\to\infty} \frac1N \left(\frac{1}{f^N_p}(\caL_N f^N_p)\right) (x)
=(e^p-1) b(x) + (e^{-p}-1) d(x)
\ee
and the corresponding Lagrangian
\be\label{lagraga}
L(x,v)=\sup_{p\in\R}\left(pv-\caH(x,p)\right)
\ee
For the trajectories of $\{X^N_t:0\leq t\leq   T\}$, we have
\be\label{largedevpoisson}
\pee (X^N_.\approx \gamma)\approx e^{-N\int_0^T L (\gamma_s, \dot{\gamma}_s)}ds
\ee
where the informal notation has to be interpreted as usual in the sense
of the large deviation principle.

The equations for the optimal trajectories, i.e. for the minimizers of the ``action''
\be\label{action}
\caI(\gamma)= i(\gamma_0) + \int_0^T L (\gamma_s,\dot{\gamma}_s) ds
\ee
can now more conveniently
be written in terms of the Hamiltonian (the Lagrangian is a more
complicated expression to deal with).

Introducing the canonical coordinates $(x,p)$
we have the Hamilton equations, together with the terminal condition
and the open-end condition corresponding to the choice
of the distribution of $X^N_0$.
\beq\label{hamiboy}
\dot{x}_t&=&\frac{\partial H}{\partial p} (x_t,p_t)= b(x)e^p-d(x)e^{-p}
\nonumber\\
\dot{p}_t &=&-\frac{\partial H}{\partial x} (x_t,p_t)=-b'(x)(e^p-1)-d'(x)(e^{-p}-1)
\eeq
with conditions
\beq
x_T&=&0
\nonumber\\
p_0 &=& i'(x_0)= 4x_0(x_0^2-a^2)
\eeq
The total ``energy'' is a constant of motion along minimizing
trajectories, so we put $H(x,p)=E$ and we can rewrite the Hamilton
equations \eqref{hamiboy}
\beq
E+b(x)+d(x)+\dot{x}&=& 2b(x) u
\nonumber\\
E+b(x)+d(x)-\dot{x}&=& 2d(x) u^{-1}
\eeq
where $u= e^{p}$.
Which leads to
\be\label{closed}
\dot{x}^2 = E^2 + 2E (b(x)+d(x)) + (b(x)-d(x))^2
\ee
So we can think now of the cost of a trajectory
as a function of two parameters: the starting
point and the energy
$(x_0,E)$. Zero-energy correspond to the ``typical trajectory''
following the limiting differential equation
$\dot{x}= b(x)-d(x)$, which means that
the cost of the Lagrangian part of the rate function is zero, and
only the cost due to the starting point $x_0$ has
to be paid. Non-zero energy
trajectories have a strictly positive cost of the
Lagrangian part of the rate function. The additional terminal condition $X_T=b$ will eliminate
one of these variables (e.g. $E$), so that we can think of the cost of
the trajectory as a function of the single variable (e.g.) $x_0$.

We now concentrate on three important particular cases.

\subsection {\bf Constant birth and death rates}
If $b$ and $d$ do not depend on $x$, then the equation for the momentum
shows that $p_t=C$, hence we have linear Euler-Lagrange trajectories, and
correspondingly the same analysis and phenomena as in the Brownian motion case
of the previous section.
\subsection{ Mean-field independent spin flip}
A special case, corresponding to independent spin flip dynamics
is
$b(x)=(1-x)$, $d(x)=(1+x)$.
Moreover, the $x$-variable is now restricted to $[-1,1]$.
As in the case $x\in\R$ we assume that initially, $x_0$ is distributed
according to a measure $\mu_n(dx)$ on $[-1,1]$ satisfying the large
deviation principle with the non-convex rate function
\eqref{rate} for $x\in [-1,1]$ and $+\infty$ otherwise. In particular,
$a\in (0,1)$.

The Hamilton equations then read
\beq
\dot{x}&=& -x(e^p+e^{-p}) + e^p-e^{-p}
\nonumber\\
\dot{p}&=&e^p-e^{-p}
\eeq
Taking the derivative w.r.t.\ time of the first equation and using
the second equation leads to elimination of $p$, and the simple second
order equation for $x$:
to
\be
\frac{d^2 x}{dt^2}= 4x
\ee
with solutions
\[
x(t)= C_1e^{2t} + C_2 e^{-2t}
\]
where $C_1, C_2$ are determined by the open-end condition and the terminal condition.
This case was treated before in the context of the Curie-Weiss model
subjected to independent spin flips in \cite{kuel}, \cite{leny}.

The equation for the momentum can be integrated and gives
\[
\tanh (p_t/2)= \pm C e^{2t}
\]
Furthermore, since
\[
E= (1-x) (e^p-1) + (1+x)( e^{-p}-1)
\]
is a constant of motion, we find as possible solutions for $x$, using that
$x_T=0$:
\[
x_t =\pm \sqrt{E/4(1+E/4)} \left(e^{2(t-T)}-e^{2(T-t)}\right)
\]
In particular, as in the Brownian motion case, the zero-energy trajectory
($E=0$) yields $x_t=0$. The relation between the energy, initial position
and initial momentum is
\[
p_0= \log\left(\frac{2+E +\sqrt{(2+E)^2 - 4(1-x_0^2)}}{2(1-x_0)}\right)
\]
Zero-energy thus corresponds to zero initial momentum and zero initial position.

In general, the initial points are symmetrically distributed around
the origin and related to the energy via
\[
x_0= \pm \sqrt{E/4(1+E/4)} \left(e^{-2T}-e^{2T}\right)
\]
Whether or not a non-zero energy solution can be the minimizer is determined
by the open-end condition:
\be\label{coribo}
p_0= i'(x_0)= 4x_0(x_0^2-a^2)
\ee
This can be viewed now as an equation for $E$. For small $T>0$,
\[
x_0= x_0(E,T)\approx C(E) T, p_0= p_0(E,T)\approx cE
\]
which implies that a non-zero energy solution of \eqref{coribo}
can not exist for
small $T$.
For large $T$, a non-zero energy solution exists
yielding two symmetrically solution for $x_0$.

Alternatively, the trajectory cost $C_T(\gamma_0)$ of a trajectory
starting at $\gamma_0$ ending up at time $T$ at $b=0$
has the following important properties
\ben
\item Symmetry: $C_T(-\gamma_0)=C_T(\gamma_0)$
\item Small time behavior: $\lim_{T\to 0} C_T(\gamma_0)=\infty$ for all $\gamma_0\not=0$
\item Large time behavior: $\lim_{T\to\infty} C_T (\gamma_0)=0$ for all $\gamma_0$
\een
From these properties it follows that for small $T$ there are no
bad points, and for large $T$ zero is the unique bad point. Notice that contrary to the
Curie Weiss model situation analyzed in \cite{kuel} there are no non-neutral
(non-zero) bad configurations due to the fact that the rate function of the initial
measure is here simply a fourth order polynomial.
\subsection{Independent spin-flips in a field}
This corresponds to the choice $b(x)= \gamma (1-x)$,
$d(x)= (1+x)$, $x\in [-1,1]$. Here $\gamma>1$ corresponds
to a bias in the plus direction (positive magnetic field).
The limiting deterministic
trajectory is given by
\[
\frac{dx_t}{dt}= -(1+\gamma) x_t + (\gamma-1)
\]
\be\label{zerengam}
x_t = x_0 e^{-(1+\gamma)t} + \frac{\gamma-1}{1+\gamma} \left(1-e^{-(1+\gamma)t}\right)
\ee
This is the zero-energy trajectory starting from $x_0$.

Using \eqref{closed} we find that for a given energy $E$, the solution for $x$ is of the form
\be\label{xfieldex}
x_t= x(E,C,t)= C_1 e^{ t(1+\gamma)} + C_2 e^{- t(1+\gamma)} + C_3
\ee
with
\beq\label{constants}
C_1 &=&
\frac{ C(2+2\gamma +E)}{(\gamma+1)^2}
\nonumber\\
C_2 &=&
\frac{E\gamma}{-C(1+\gamma)^2)}
\nonumber\\
C_3 &=& \left(\frac{(E+1+\gamma)(\gamma-1)}{(1+\gamma)^2}\right)
\eeq
where $C$ is an integration constant.
\br
\ben
\item Remark that for $E=0$ $C_3= (\gamma-1)(1+\gamma)^{-1}$ which corresponds to the limiting value
of the zero energy trajectory.
\item If $\gamma=1$, and $E\not=0$ we find $C_3=0$ and recover the solution of the form
$C_1e^{2t}+ C_2e^{-2t}$ corresponding to the optimal trajectories of the independent spin flip dynamics.
\een
\er
The general form of an optimal trajectory arriving at time $T$ at $x_T=b$ and
starting from $x_0=\gamma_0$
is
\[
x(t)= (b-C_3) \frac{\sinh (\delta t)}{\sinh (\delta T)}
+ (\gamma_0-C_3)\frac{\sinh (\delta (T-t))}{\sinh (\delta T)} +C_3
\]
with $\delta= (1+\gamma)$ and
where $C_3$ is given in \eqref{constants}. Notice the analogy with the case of the Ornstein-Uhlenbeck process
in a constant field \eqref{badorn}.
As in that case, the bad point is time-dependent and given by
\[
b= \frac{\gamma-1}{\gamma+1} \left(1-e^{-\delta T}\right)
\]
which is the point at which the limiting deterministic dynamics arrives at time $T$ when started
from $x_0=0$.
The trajectory cost $C_T(\gamma_0)$ to arrive at this bad point satisfies the same properties as the trajectory cost
$C_T(\gamma_0)$ of the previous subsection (zero field case).
Hence, for $T$ large two minimizing $\gamma_0$ of the total cost function
appear which correspond to two optimal trajectories.

\section{ Acknowledgement}
We thank Aernout van Enter and Olaf de Leeuw for useful discussions and suggestions.

\end{document}